% plain TeX file of the paper ARGUMENT PRINCIPLE AND HOLOMORPHIC EXTENDIBILITY
% written by Josip Globevnik
% prepared on April 25, 2004
% submitted for publication in Journal d'Analyse Mathematique
% accepted on May 21, 2004
\magnification 1200
\def\R{{\rm I\kern-0.2em R\kern0.2em \kern-0.2em}}
\def\N{{\rm I\kern-0.2em N\kern0.2em \kern-0.2em}}
\def\P{{\rm I\kern-0.2em P\kern0.2em \kern-0.2em}}
\def\B{{\rm I\kern-0.2em B\kern0.2em \kern-0.2em}}
\def\C{{\bf \rm C}\kern-.4em {\vrule height1.4ex width.08em depth-.04ex}\;}

\def\D{{\Delta}}

\def\z{{\zeta}}

\def\cH{{\cal H}}
\def\A{{A(D)}}
\def\e{{e^{i\omega}}}
\def\G{{\Gamma}}
\def\alpha{{\beta}}
\font\ninerm=cmr8
\ 
\vskip 6mm
\centerline {\bf THE ARGUMENT PRINCIPLE AND HOLOMORPHIC EXTENDIBILITY}
\vskip 4mm
\centerline{Josip Globevnik}
\vskip 4mm
{\noindent \ninerm ABSTRACT\ \  Let $D$ be a bounded domain in the complex plane 
whose boundary consists of finitely many 
 pairwise disjoint simple closed curves. Give $bD$ the standard orientation and let 
$A(D)$ be the algebra of all continuous 
functions on $\overline D$ which are holomorphic on $D$. In the paper we prove that a 
continuous function $f$ on $bD$ extends to a function in $A(D)$ if and only if for each 
$g\in A(D)$ such that $f+g\not=0 $ on $bD$ the change of argument 
of $f+g$ along $bD $ is nonnegative. } 
\vskip 6mm
\bf 1.\ Introduction and the main result \rm 
\vskip 2mm
H.\ Alexander and J.\ Wermer [AW] obtained a characterization of boundaries 
of a\-na\-lytic varieties in terms of a generalized argument principle. Their results 
brought new results also into the classical function theory of one variable. For 
instance, E.\ L.\ Stout [S1] observed that their results imply that a smooth function 
on the unit circle 
$b\D $ extends holomorphically through the open unit disc $\D$ if and only 
if for each polynomial $Q$ of two variables such that 
$Q(z,f(z))\not= 0\ (z\in b\D )$ the change of argument of $z\mapsto Q(z,f(z))$ 
around $b\D $ is nonnegative, and proved that this holds for continuous functions. 
J.\ Wermer [W] showed that for smooth functions on $b\D $ there is a better result: 
a smooth function $f$ on $b\D $ extends holomorphically 
through $\D $ if and only if for each polynomial $P$ such that $f+P\not= 0$ on $b\D $ 
the change of argument of $f+P$ around $b\D $ is nonnegative. For continuous functions
this was proved by the author in [Gl] so that we have the following 
characterization of the disc algebra in terms of the 
argument principle:
\vskip 2mm \noindent\bf Theorem 1.0\ \rm [Gl]\ \it A continuous function $f$ on $b\D $ 
extends holomorphically through $\D $ if and only if for each polynomial $P$ 
such that $f+P\not= 0$ on $b\D $ 
the change of argument of $f+P$ around $b\D $ is nonnegative. 
\vskip 2mm
\noindent \rm In the present paper we prove that the analogous theorem holds for 
multiply connected domains.

Let $D\subset\C$ be a bounded domain whose boundary consists of
finitely many pairwise disjoint simple 
closed curves. We give $bD$ the standard orientation. Denote by $A(D)$ the algebra of all 
continuous functions on $\overline D$ which are holomorphic on $D$. Our main result is
\vskip 2mm
\noindent\bf THEOREM 1.1\ \it A continuous function $f$ on $bD$ extends to a function 
in $\A$ if and only if for each $g\in\A$ such that $f+g\not=0$ on $bD$ the change 
of argument of $f+g$ along $bD$ is nonnegative. \rm
\vskip 2mm
If the condition in Theorem 1.1 holds for all $g$ belonging to a dense subset of $\A$ then it 
holds for all $g\in\A$. Thus, since rational functions with poles outside 
$\overline D$ are dense in $A(D)$ [S2, 23], it is enough to assume that the condition in 
Theorem 1.1  
holds for rational functions with poles outside $\overline D$. 
The only if part of the theorem is an obvious consequence of the 
argument principle. In fact, if $f$ admits an extension $\tilde f\in A(D)$ then the change of 
 argument of $f+g$ along $bD$ equals $2\pi $ times the number of zeros of $\tilde f+g$ 
in $D$. 

\vskip 4mm
\bf 2.\ Preliminaries \rm
\vskip 2mm
Every bounded domain $D\subset \C$ whose boundary consists of finitely 
many pairwise disjoint simple closed 
curves is biholomorphically equivalent to a domain $D^\prime$ bounded by finitely many 
pairwise disjoint circles [Go]. Moreover, every biholomorphic map 
$\Phi\colon D\rightarrow D^\prime $ extends to 
a homeomorphism $\tilde\Phi\colon \overline D\rightarrow\overline{D^\prime}$ 
(see the proof in [CL,\ pp.\ 46-49] which works 
also for multiply connected domains, bounded by finitely many pairwise disjoint 
simple closed curves). 
Thus, with no loss 
of generality assume that $D$ is bounded by finitely many pairwise disjoint circles.     

In general, not every real-valued harmonic function $u$ on $D$ 
is the real part of a holomorphic function on $D$. If there is a harmonic function $v$ on $D$ 
such that $u+iv$ is holomorphic on $D$ then we call $v$ a \it conjugate of $u$. \rm If $D$ 
is simply connected then every harmonic function on $D$ has a conjugate on $D$. Let $f$ be a 
complex valued harmonic function on $D$. Write $f= p+iq$ with $p,\ q$ real. We will say that 
$f$ has a conjugate on $D$ if $p$ has a conjugate $r$ on $D$ and $q$ has a conjugate $s$ on 
$D$, and we will call the function $r-is$ a conjugate of $f$ on $D$.  This 
happens if and only if $f=F+\overline G$ where F and G are holomorphic functions on $D$. In 
fact, if $P=p+ir$ and $Q=q+is $ then $P$ and $Q$ are holomorphic functions on $D$ and  
$F=(P+iQ)/2$ and $G=(P-iQ)/2$. 

Given a continuous function $\Phi $ on $bD$ there is a continuous extension of $\Phi$ to 
$\overline D$ which is harmonic on $D$ and which we will denote by $\cH (\Phi )$; moreover, if 
$\Phi\in C^\infty (bD)$ then $\cH (\Phi )\in C^\infty (\overline D)$ [B, p.53]. 
$\cH $ is a linear 
map from $C(bD)$ to the space of continuous functions on $\overline D$ which
are harmonic on $D$. If D is simply 
connected then $\cH (\Phi )$ has a conjugate harmonic function on $D$ which is also 
in $C^\infty (\overline D)$ provided  that $\Phi\in C^\infty (bD)$ [B, p.91]. If $\Omega _1$ 
and $\Omega _2$ are bounded, simply connected domains with boundaries of class $C^\infty $ 
then a biholomorphic map $\Phi$ from $ \Omega _1 $ to $ \Omega _2$ extends to a smooth map 
from $\overline{\Omega _1}$ to $\overline{\Omega _2}$ [B, p.28]. Applying this locally along $bD$ 
and using 
the preceding discussion we see that if $f$ is a harmonic function on $D$ that has a 
conjugate on D then the conjugate extends smoothly to $\overline D$ provided that $f$ extends 
smoothly to $\overline D$. We summarize this in
\vskip 2mm
\noindent
\bf LEMMA 2.1\ \it If $f\in C^\infty (bD)$ is such that $\cH (f)$ has a conjugate on 
$D$ then both $\cH (f)$ and its conjugate extend smoothly to $\overline D$. \rm 
\vskip 2mm
Every 
harmonic 
function $f$ on $D$ is real-analytic on $D$ so, if $f$ is holomorphic on an open 
subset of $D$
then it is holomorphic on $D$. We will need the following fact which can be found 
as an exercise in [R]. 
\vskip 2mm
\noindent \bf LEMMA 2.2\ \it Let $f$ be a harmonic function on $D$ such that $z\mapsto 
zf(z) $ is harmonic on a nonempty open set $U\subset D$. Then $f$ is holomorphic on $D$. \rm
\vskip 1mm
\noindent \bf Proof.\ \rm By the preceding discussion it is enough to prove that $f$ is holomorphic on 
a disc $\Omega \subset D$. Since $f$ is harmonic on $D$ there is a disc $\Omega\subset U$ such that 
$f=P+\overline Q$ on $\Omega $ where $P$ and $Q$ are holomorphic functions on $\Omega $. By our assumption, the 
function $z\mapsto zf(z) = zP(z)+z\overline{Q(z)}$ is harmonic on $\Omega $ so
$$
{{\partial^2}\over{\partial z\partial\overline z}}[zP(z)+z\overline{Q(z)}] =0\ \ (z\in\Omega)
$$
which implies that $\overline{Q^\prime (z)}=0\ (z\in\Omega)$. Thus, $Q$ is constant on $\Omega $ 
and consequently
$f$ is holomorphic on $\Omega $. This completes the proof.
\vskip 4mm
\bf 3.\ A new proof in the case of a disc \rm
\vskip 2mm
The proof in [Gl] does not 
generalize to multiply connected domains. In this section we give a new, different proof 
of the theorem in the case when 
$D$ is a disc which we later generalize to multiply connected domains. 

Write $\D = \{\z\in\C\colon |\z|<1\}$. Throughout this section, $D=\D$, Denote by $Z$ the 
identity function:\ \ $Z(z)=z \ 
(z\in bD)$. Assume that $f\in C(b\D)$  does not extend to a function 
from the disc algebra $A(\D)$. Then there is an $a\in\D$ such that 
$$
\cH[(Z-a)f](a)\not= 0.
\eqno (3.1)
$$
To see this, suppose for a moment that $\cH[(Z-a)f](a)= 0 \ (a\in\D)$. This implies that 
$\cH (Zf)(a)=a\cH (f)(a)\ (a\in\D)$. 
In particular, the function $a\mapsto a\cH (f)(a)$ is harmonic on $\D$. Since $\cH (f)$ is harmonic on $D$ 
Lemma 2.2 implies that $\cH (f)$ is holomorphic on $\D$ so $f$ extends holomorphically through $\D$,
a contradiction. This proves that there is an $a\in\D$ such that (3.1) holds. 

With no loss of generality, replacing $f$ with $\e f$, $\omega\in\R$, if necessary, 
we may assume that 
$$
\Re\{\cH[(Z-a)f](a)\} =\alpha\not= 0.
\eqno (3.2)
$$
For easier understanding we complete the proof first under the additional 
assumption that $f $ is smooth. 
Suppose for a moment that $f$ is smooth. The function $z\mapsto \Re \{\cH[(Z-a)f](z)\}
-\alpha $ 
is continuous on $\overline\D$, harmonic on $\D $ and has smooth boundary values 
$\Re [(z-a)f(z)]-\alpha \ (z\in b\D)$. Hence by Lemma 2.1 
there is a function $g\in A(\D )$ such that 
$$
\Re [g(z)] = \Re\{\cH[(Z-a)f](z)\}-\alpha  \ \ (z\in\overline\D ).
$$
By (3.2) we have $\Re g(a)=\Re\{\cH[(Z-a)f](a)\}-\alpha = 0$ so by adding an imaginary constant to $g$ if necessary we may assume that $g(0)=0$ so 
$g(z)=(z-a)h(z)\ (z\in\overline\D )$ where $h\in A(\D )$. Consider the function 
$z\mapsto G(z) = (z-a)f(z)-g(z)\ (z\in b\D)$. We have $\Re [G(z)]=\alpha \ (z\in b\D)$ 
which, since $\alpha\not=0$, implies that 
$G\not= 0$ on $b\D$ and that the change of argument of $G$ 
around $b\D $ equals zero. Since $a\in\D$ the change of argument of $z\mapsto (z-a)$ around $b\D $ equals 
$2\pi $. Since 
$G(z)= (z-a)[f(z)-h(z)]\ (z\in b\D)$ it follows that $f-h\not=0$  on  $b\D$ and that 
the change of argument of 
$f-h$ around $b\D $ is negative. Since $h\in A(\D )$ this completes the proof 
in the special case when $f$ is smooth. 

In general, we have to approximate $f$ by smooth functions as follows:

\noindent Let $f_1$ be a smooth function on $b\D $ such that 
$$
|(z-a)[f(z)-f_1(z)]|<|\alpha |/4\ \ (z\in b\D ).
\eqno (3.3)
$$
Write $\alpha _1= \Re\{\cH[(Z-a)f_1](a)\}$. By (3.3) and by the
maximum principle for the real harmonic function $\Re \{\cH[(Z-a)(f_1 - f)]\}$ we have 
$$
|\alpha _1-\alpha |<|\alpha|/4 .
\eqno (3.4)
$$
The function $z\mapsto \Re \{\cH [(Z-a)f_1](z)\}-\alpha_1$ is continuous on $\overline D$, 
harmonic on $\D$ and has smooth boundary values $\Re [(z-a)f_1(z)]-\alpha _1\ (z\in b\D )$. 
Hence by Lemma 2.1 there is a function $g_1\in A(\D )$ such that 
$$
\Re [g_1(z)]= \Re \{\cH[(Z-a)f_1](z)\}-\alpha _1\ (z\in\overline\D ).
\eqno (3.5)
$$
Clearly $\Re [g_1(a)] = \Re \{\cH[(Z-a)f_1](a)\}-\alpha _1= 0$ so by adding an imaginary constant to $g_1$ if 
necessary we may assume that $g_1(a)=0$ so that $g_1(z)=(z-a)h(z)\ (z\in\overline\D)$ where 
$h\in A(\D )$. 

Consider the function 
$
z\mapsto G(z)= (z-a)f(z)-g_1(z)\ (z\in b\D )
$.  
By (3.5) we have $\Re G(z) = \alpha +\Re [(z-a)
(f(z)-f_1(z))] + (\alpha_1 -\alpha) \ (z\in b\D )$ which, by (3.3) and (3.4) implies that 
$$
|\Re G(z)-\alpha |< |\alpha |/4+|\alpha |/4 = |\alpha |/2\ \ (z\in b\D )
$$
so 
$$
\alpha - |\alpha |/2<\Re G(z)<\alpha + |\alpha |/2   \ (z\in \D ),
$$
which, since $\alpha\not=0$, implies that $G\not=0$ on $b\D $ and that the change of 
argument of $G$ around $b\D $ is zero. We now repeat the reasoning 
from the proof in the smooth case to conclude that $f-h\not=0 $ on $b\D $ and that 
the change of argument of $f-h$ around $b\D $ is negative. 
Since $h\in A(\D )$ this completes the proof.
\vskip 4mm
\bf 4.\ Harmonic functions and their conjugates \rm
\vskip 2mm
The main problem in generalizing the proof in Section 3 to multiply connected 
domains $D$ is that 
in general, a harmonic function on $D$ has no conjugate on $D$.

We have assumed that  $D\subset\subset \C$ is a domain bounded by 
pairwise disjoint circles. Denote these circles by  
$\G _1,\G _2,\cdots ,\G_n$ where $\G _n$ is the boundary of the unbounded 
component of $\C\setminus\overline D$. 
For each $k, 1\leq k\leq n$, the \it harmonic measure function \rm $\omega _k$ 
is the continuous function on  $\overline D$, harmonic on $D$ which satisfies 
$\omega _k\equiv 1$ on $\G _k$ 
and 
$\omega _k \equiv 0 $ on  $\G _j,\ 1\leq j\leq n, j\not=k$. By the preceding discussion 
each $\omega _k, 1\leq k\leq n$, is smooth on $\overline D$. We have  $\sum_{k=1}^n\omega _k
\equiv 1$ on $\overline D$.

For each $k,\ 1\leq k\leq n-1$, let $\gamma _k$ be a circle with the same center as $\G _k$ and with a a slightly 
larger radius, and let $\gamma_n $ be a circle with the same center as $\G _n$ and 
with a slightly smaller radius so that the circles $\gamma _k, \ 1\leq k\leq n$ bound a domain $D^\prime $, slightly 
smaller than $D$, whose closure is contained in $D$. We give each $\gamma _k$ the orientation induced by 
the standard orientation of $bD^\prime $. 

Let $u$ be a real-valued harmonic function on $D$. 
A conjugate $v$ of $u$ has to satisfy the Cauchy-Riemann equations
$$
{{\partial v}\over{\partial x}}= -{{\partial u}\over{\partial y}}\ \ \  \ \ \ \ 
{{\partial v}\over{\partial y}}= {{\partial u}\over{\partial x}} .\ \ \ 
$$
This system is always solvable for $v$ locally. It is solvable for $v$ on $D$ if and only if 
$$
\int _{\gamma _k}\Bigl( 
-{{\partial u}\over{\partial y}} dx + 
{{\partial u}\over{\partial x}} dy\Bigr) = 0\ \ \ \ (1\leq k\leq n-1).
\eqno (4.1)
$$
Since 
$$
2\int _{\gamma _k}{{\partial u}\over{\partial z}} dz = i \int _{\gamma _k}\Bigl( 
-{{\partial u}\over{\partial y}} dx + 
{{\partial u}\over{\partial x}} dy\Bigr) \ \ \ (1\leq k\leq n-1)
$$
(4.1) holds if and only if 
$$
\int _{\gamma _k}{{\partial u}\over{\partial z}} dz = 0\ \ (1\leq k\leq n-1).
\eqno (4.2)
$$
If $u$ is a real valued harmonic function on $D$ then there are real constants
$c_1, c_2,\cdots ,c_{n-1}$ such that 
the harmonic function $u+\sum_{j=1}^{n-1}c_j\omega _j$ has a conjugate on $D$. For 
this to happen we must have 
$$
\int _{\gamma _k} {{\partial }\over{\partial z}}\Bigl[ u+\sum_{j=1}^{n-1}c_j\omega _j \Bigr] dz = 0
\ \ (1\leq k\leq n-1),
$$
that is, 
$$
\sum _{j=1}^{n-1} c_j \int_{\gamma _k}{{\partial \omega_j}\over{\partial z}} dz = 
-\int _{\gamma _k}{{\partial u}\over{\partial z}} dz 
\eqno (4.3)
$$
The system (4.3) has  a unique solution since the matrix
$$
\Bigl[\int_{\gamma _k}{{\partial\omega_j}\over{\partial z}} dz\Bigr]_{1\leq j,k\leq n-1} =
\Bigl[\int_{\G _k}{{\partial\omega_j}\over{\partial z}}dz \Bigr]_{1\leq j,k\leq n-1} 
$$
is known to be nonsingular [B, p.82]. (In the last equality we used the fact that 
all functions 
${{\partial\omega_j}\over{\partial z}}$ are smooth on $\overline D$ and holomorphic on $D$.) 
The Poisson formula implies that 
given $\varepsilon >0$ and a compact set $K\subset D$ there is a $\delta >0$ such that 
$\vert{{\partial u}\over{\partial z}}\vert<\varepsilon $ on $K$ whenever $u$ is a real harmonic 
function on $D$ such that $|u|<\delta $ on $D$. Thus, given $\varepsilon >0$ 
there is a $\delta >0$ such that 
$$
\Bigl\vert \int _{\gamma _k}{{\partial u_1}\over{\partial z}}dz -\int_{\gamma _k}
{{\partial u}\over{\partial z}}dz \Bigr\vert <\varepsilon \ \ (1\leq k\leq n-1)
$$
provided that $u_1$ is a real harmonic function on $D$ satisfying  
$|u_1(z)-u(z)|<\delta \ (z\in D)$. 

The preceding discussion gives
\vskip 2mm
\noindent\bf LEMMA 4.1\ \it Given a harmonic function $f$ on $D$ there is a unique (n-1)-tuple 
$c_1(f),$ $\cdots , $ $c_{n-1}(f)$ of complex numbers such that  $f+\sum_{j=1}^{n-1} c_j(f)\omega _j$
has 
a conjugate on $D$. These numbers depend continuously on $f$ in the sup norm on $D$.  
\vskip 4mm
\bf 5.\ Proof of Theorem 1.1, Part 1 \rm
\vskip 2mm
Let $D$ be a domain bounded by pairwise disjoint circles $\G _1,\cdots ,\G _n$ 
where $\G _n$ is the boundary of the unbounded component of $\C\setminus\overline D$. Section 3 contains 
the proof of Theorem 1.1 in the case $n=1$ so assume that $n\geq 2$. 

Let $f$ be a continuous function on $bD$ which does not extend holomorphically through $D$. Define 
$$
A(a,f)=\cH [(Z-a)f](a) = \cH (Zf)(a)-a\cH (f)(a)\ \ (a\in \overline D).
$$
Since $\cH (f)$ is not holomorphic Lemma 2.2 implies that
$$
\{a\in D\colon\ A(a,f)=0\}\hbox{\ \ is a closed, nowhere dense subset of \ }D.
\eqno (5.1)
$$
There are constants $c_k(f),\ 1\leq k\leq n-1,$ such that 
$$
\cH (f) (z) + \sum_{k=1}^{n-1} c_k(f)\omega _k(z) = P_f(z)+\overline{Q_f(z)} \ \ (z\in D)
\eqno (5.2)
$$
where $P_f$ and $Q_f$ are holomorphic functions on $D$.  Similarly, there are constants 
$d_k(f),\ 1\leq k\leq n-1,$ such that 
$$
\cH (Zf)(z)+\sum_{k=1}^{n-1} d_k(f)\omega _k(z) = R_f(z)+\overline{S_f(z)} \ \ (z\in D)
\eqno (5.3)
$$
where $R_f$ and $S_f$ are holomorphic functions on $D$. We know that the constants 
$c_k(f)$ and $d_k(f),\ 1\leq k\leq n-1$, are determined uniquely and, by the maximum 
principle for harmonic functions depend continuously on $f\in C(bD)$. We have
$$
c_k(\e f)=\e c_k(f),\ \ d_k(\e f)=\e d_k(f)\ \ (1\leq k\leq n-1,\ \omega\in \R)
\eqno (5.4)
$$
and
$$
A(a,\e f) = \e A(a,f)\ \ (\omega\in\R , \ a\in D).
\eqno (5.5)
$$
Define
$$\Phi_{a,f}(z) =\sum_{j=1}^{n-1}[d_j(f)-ac_j(f)].[\omega _j(z)-\omega _j(a)]-A(a,f)
\ \ (z\in D).
$$
For each $a\in D$ the function $\Phi_{a,f}$ is smooth on $\overline D$ and harmonic on $D$. 
By the preceding discussion for each $a\in D $ the harmonic function
$$
z\mapsto \cH[(Z-a)f](z) + \Phi _{a,f}(z)\ \ (z\in D)
$$
vanishes at $a$ and has a conjugate on $D$, that is, it is of the form $F_{a,f}+\overline{G_{a,f}}$  
where  $F_{a,f}$ and $G_{a,f}$ are holomorphic on $D$.
\vskip 4mm
\noindent\bf 6.\ The function $\Phi_{a,f}$ \rm 
\vskip 2mm
Note that for each $a\in D$ the function $\Phi_{a,f}$ is smooth on $\overline D$, 
harmonic on $D$ and constant on each component $\G _j,\ 1\leq j\leq n$, of $bD$. 
\vskip 2mm
\noindent \bf LEMMA 6.1\ \it There is an $a\in D$ such 
that $\Phi _{a,f}(z)\not=0 \ (z\in bD)$. \rm 
\vskip 1mm
\noindent\bf Proof.\ \rm Recall that for each $a\in D$ the function $\Phi _{a,f}|\G _k$ 
is constant for each $k,\ 1\leq k\leq n$; we have to prove that for some 
$a\in D$ these constants are all different from $0$. We shall prove that 
$$
\left.\eqalign{&
\hbox{for each\ }k,\ 1\leq k\leq n,\hbox{\ the set\ }\{a\in D\colon \ \Phi_{a,f}|\G _k = 0\}\cr
&\hbox{is a closed subset of\ }D\hbox{\ with empty interior.}\cr}\right\}
\eqno (6.1)
$$
Assume that we have done this. Then $\cup_{k=1}^n\{a\in D\colon \ \Phi_{a,f}|\G _k = 0\} $ is 
a closed subset of $D$ with empty interior which implies that there is an open dense subset 
of $D$ of those $a$ for which $\Phi _{a,f}(z)\not=0 \ (z\in bD)$ which 
will complete the proof. It remains to prove (6.1).

Let $1\leq k\leq n-1$. On $\G _k$ the function $\Phi _{a,f}$ is  
equal to the constant $-A(a,f)+
\sum_{j=1, j\not= k}^{n-1}[d_j(f)-ac_j(f)].[-\omega _j(a)] +
[d_k(f)-ac_k(f)].[1-\omega _k(a)] $. Since $a\mapsto A(a,f)$ is continuous on $D$ 
it follows that 
$\{ a\in D\colon\ \Phi_{a,f}|\G _k =0\}$ is a closed subset of $D$. Suppose that it contains a disc $U$. 
Then 
$$
A(a,f)=\sum_{j=1, j\not= k}^{n-1}[d_j(f)-ac_j(f)].[-\omega _j(a)] +
[d_k(f)-ac_k(f)].[1-\omega _k(a)] 
\eqno (6.2)
$$
for all $a\in U$. Since both sides of (6.2) are real-analytic in $a$ on $D$ 
it follows that (6.2) holds 
for all $a\in D$. Since both sides of (6.2) are continuous in $a$ on $\overline D$ it 
follows that (6.2) holds for all $a\in bD$. However, $A(a,f)=0\ (a\in bD)$ so 
$$
\sum_{j=1, j\not= k}^{n-1}[d_j(f)-ac_j(f)].[\omega _j(a)] =
[d_k(f)-ac_k(f)].[1-\omega _k(a)] \ \ (a\in bD).
\eqno (6.3)
$$
If $a\in \G _j,\ 1\leq j\leq n-1, j\not = k$, then $\omega _j(a)=1$ and 
$\omega _i(a)=0$ for 
all $i, \ 1\leq i\leq n,\ i\not= j$,  so (6.3) implies that 
$$
d_j(f)-ac_j(f)=d_k(f)-ac_k(f)\ \ (a\in\G _j, 1\leq j\leq n-1, j\not=k).
\eqno (6.4)
$$
If $a\in\G _n$ then $\omega _j(a) = 0$ for all
$j,\ 1\leq j\leq n-1$, including $k$, so (6.3) gives
$$
d_k(f)-ac_k(f)=0\ \ (a\in\G _n).
\eqno (6.5)
$$
Now, (6.5) implies that $d_k(f)=c_k(f)=0$ which, by (6.4) 
gives $d_j(f)=c_j(f)=0\ (1\leq j\leq n-1,\ j\not= k)$ so, by (6.2) it follows that 
$A(a,f)=0$ for every $ a\in D$ which contradicts (5.1). This proves that $\{ a\in D\colon\ \Phi_{a,f}|\G _k =0\}$ 
has empty interior for each $k,\  1\leq k\leq n-1.$

Let $k=n$. We have
$$
\Phi_{a,f}|\G _n= -A(a,f)+\sum _{j=1}^{n-1}[d_j(f)-ac_j(f)].[-\omega _j(a)]
$$ 
As before, the continuity of $a\mapsto A(a,f)$ on $D$ implies that the 
set $\{ a\in D\colon\ \Phi _{a,f}|\G _n =0\} $ is closed. Suppose that it has nonempty interior. 
As before, we get 
$$
0=A(a,f)=\sum _{j=1}^{n-1}[d_j(f)-ac_j(f)].[-\omega _j(a)] \ \ (a\in bD)
\eqno (6.6)
$$
It follows that $d_j(f)-ac_j(f)=0\ \ 
(a\in\G _j,\ 1\leq j\leq n-1\} $  which implies that 
$d_j(f)=c_j(f)=0\ \ (1\leq j\leq n-1)$ so again $A(a,f)\equiv 0\ (a\in D)$ which contradicts (5.1). 
This completes the proof. 
\vskip 4mm
\bf 7.\ Proof Theorem 1.1, Part 2\rm 
\vskip 2mm
By Lemma 6.1 there is an $a\in D$ such that the constants 
$\Phi_{a,f}|\G _k,\ 1\leq k\leq n,$ are all different from $0$. By (5.4) and (5.5) we have
$\Phi_{a,\e f}=\e \Phi _{a,f}\ (\omega\in\R)$ so replacing $f$ by $\e f$ if necessary 
we may assume with no loss of generality that
$$
(\Re \Phi_{a,f})|\G _k = \alpha_k\not= 0\ \ (1\leq k\leq n).
\eqno(7.1)
$$
To make the proof easier to understand we first complete it under the 
assumption that $f$ is smooth. Assume that $f$ is smooth. In this case, by Lemma 2.1, 
$$
\cH [(Z-a)f](z)+\Phi _{a,f}(z) = F_{a,f}+
\overline {G_{a,f}(z)}\ \ (z\in \overline D)
$$ 
where  $F_{a,f}$ and $G_{a,f}$ belong to $A(D)$ so
$$
\Re \{\cH[(Z-a)f](z)+\Phi _{a,f}(z)\} = \Re [g(z)]\ \ (z\in\overline D)
$$
where $g=(F_{a,f}+G_{a,f})/2 \in A(D)$. Clearly $\Re [g(a)]=0$ so by 
adding an imaginary constant to $g$ if necessary we may assume that $g(0)=0$ so 
$g(z)=(z-a)h(z)\ (z\in\overline D)$ where $h\in A(D)$. 

Consider the function $z\mapsto G(z) = (z-a)f(z)-g(z)\ (z\in bD)$. We have 
$\Re [G(z)]= -\Re \Phi_{a,f}(z)\ \ (z\in bD)$. By (7.1) for 
each $k,\ 1\leq k\leq n,$ the expression on the right is a
nonzero constant on $\G _k$ which implies that for each
$k,\ 1\leq k\leq n$, the change of argument of $z\mapsto G(z)=(z-a)[f(z)-h(z)]$ along 
$\Gamma _k$ equals $0$. Since $a\in D$ the change of argument of $z\mapsto (z-a)$
along each $\G _k,\ 1\leq k\leq n-1$, is zero, and the change of argument 
of $z\mapsto (z-a)$ along $\G _n$ is $2\pi$. Thus, the change of argument of $z\mapsto 
f(z)-h(z)$ along $\G _n$ equals $-2\pi $. So, the change of argument 
of $z\mapsto f(z)-h(z)$ along $bD$ is negative. Since $h\in A(D)$ this completes the proof in the 
case when $f$ is smooth.

In the case of general $f$ we have to approximate $f$ by smooth functions. 
We already know that the constants $c_k(f)$ and $d_k(f)$ depend 
continuously on $f\in C(bD)$. Further, for our fixed $a\in D$ the maximum 
priciple for harmonic functions implies that $A(a,f)= \cH [(Z-a)f](a)$ also depends 
continuously on $f\in C(bD)$. It follows that $\Phi _{a,f_1}$ is uniformly arbitrarily close to 
$\Phi _{a,f}$ on $\overline D$ provided that $f_1\in C(bB)$ is sufficiently close to $f$. Fix $\varepsilon$, 
$$
0<\varepsilon< (1/4) \min \{|\alpha _k|\colon\ 1\leq k\leq n\} 
\eqno (7.3)
$$
and let $f_1$ be a smooth function on $bD$ which is so close to $f$ that
$$
|\Phi _{a,f_1}(z) - \Phi _{a,f}(z)|<\varepsilon \ \ (z\in\overline D)
\eqno (7.4)
$$
and
$$
|(z-a)[f_1(z)-f(z)]|<\varepsilon \ \ (z\in bD).
\eqno (7.5)
$$
As before, since $f_1$ is smooth, Lemma 2.1 applies to show that  
$$
\Re \{\cH [(Z-a)f_1](z)+\Phi _{a,f_1}(z)\} = \Re [g(z)]\ \ (z\in\overline D)
\eqno (7.6)
$$ 
where $g\in A(D)$ satisfies $g(a)=0$. Consider the 
function $z\mapsto G(z)=(z-a)f(z)-g(z)$. We have  
$\Re [G(z)] = \Re \{ (z-a)[f(z)-f_1(z)]\} + \Re [\Phi_{a,f}(z)] +
\Re [\Phi_{a,f_1}(z)-\Phi_{a,f}(z)]$ so by (7.4) and (7.5) it follows that 
$$
\bigl\vert \Re [G(z)] - \Re [\Phi _{a,f}(z)]\bigr\vert < 2\varepsilon\ \ (z\in bD)
$$
so by (7.1) it follows that
$$
\bigl\vert \Re [G(z)] - \alpha _k\bigr\vert <2\varepsilon \ \ (z\in \G _k,\ 1\leq k\leq n),
$$
which, by (7.3), implies that 
$$
\alpha _k - |\alpha _k|/2 < \Re [G(z)]<\alpha _k + |\alpha _k|/2\ \ 
(z\in \G _k,\ 1\leq k\leq n).
$$
Since $\alpha _k\not=0\ (1\leq k\leq n)$, it follows that for each $k, \ 1\leq k\leq n$, 
the change of argument of $G$ along $\G _k$ is zero. 
Now we conclude the proof as in the smooth case.
\vskip 5mm
\noindent \bf Acknowledgement\ \rm  A part of the work whose results are presented 
here was done during the author's visit at the University of Oslo in November 2003. 
The author 
is grateful to Erik L\o w and Niels \O vrelid for their hospitality. 

This work was supported 
in part by the Ministry of Education, Science and Sport of Slovenia 
through research program Analysis and Geometry, Contract No.\ P1-0291. 
\vfill
\eject
\centerline{\bf REFERENCES} 
\vskip 5mm
\noindent [AW]\ H.\ Alexander and J.\ Wermer: Linking 
numbers and boundaries of varieties. 

\noindent Ann.\ Math.\ 151 (2000) 125-150
\vskip 2mm
\noindent [B]\ \ S.\ Bell: \it The Cauchy transform, Potential Theory, and 
Conformal Mapping. \rm CRC Press, 
Boca Raton, 1992
\vskip 2mm
\noindent [CL]\ \ E.\ F.\ Collingwood and A.J.Lohwater: \it The Theory of Cluster Sets. \rm

\noindent Cambridge University Press, Cambridge 1966
\vskip 2mm
\noindent [Gl]\ \ J.\ Globevnik:\ Holomorphic extendibility and the argument principle.

\noindent To appear in "Complex Analysis and Dynamical Systems II (Proceedings
  of a conference held in honor of Professor Lawrence Zalcman's sixtieth
  birthday in Nahariya, Israel, June 9-12, 2003)", Contemp.\ Math. 
  [http://arxiv.org/abs/math.CV/0403446]

\vskip 2mm
\noindent [Go]\ \ G.\ M.\ Goluzin: \it Geometrische Funktionentheorie. \rm
\rm VEB Deutscher Verlag 
der Wissenschaften, Berlin 1957
\vskip 2mm
\noindent [R]\ W.\ Rudin:\ \it Real and Complex Analysis. \rm 
McGraw-Hill, New York, 1970
\vskip 2mm
\noindent [S1]\ E.\ L.\ Stout:\ Boundary values and mapping degree.

\noindent Michig.\ Math.\ J.\ 47 (2000) 353-368
\vskip 2mm
\noindent [S2]\ E.\ L.\ Stout: \it The Theory of Uniform Algebras.\ \rm Bogden and Quigley, 
Tarrytown -on-Hudson, N.Y. 1971
\vskip 2mm
\noindent [W] J.\ Wermer: The argument principle and boundaries of analytic varieties.

\noindent Oper. Theory Adv. Appl., 127, Birkhauser, Basel, 2001, 639-659
\vskip 20mm

\noindent 
Institute of Mathematics, Physics and Mechanics

\noindent 
University of Ljubljana

\noindent Ljubljana, Slovenia

\noindent josip.globevnik@fmf.uni-lj.si

\bye